\documentclass[a4paper,11pt,oneside]{article}
\pdfoutput=1 
\usepackage[english]{babel}
\usepackage[dvips]{graphicx}
\usepackage[svgnames]{xcolor}
 
\graphicspath{{./Figures/}}

\usepackage{amssymb}
\usepackage[tbtags]{amsmath}
\usepackage{amsthm}
\usepackage{icomma}

\theoremstyle{definition}
\newtheorem{definition}{Definition}

\theoremstyle{plain}

\usepackage[T1]{fontenc}
\usepackage[osf,sc]{mathpazo}
\usepackage[euler-digits,euler-hat-accent]{eulervm}

\usepackage[top=32mm]{geometry}

\usepackage{caption}
\usepackage{booktabs}
\usepackage[flushleft]{threeparttable}

\usepackage{tabularx}

\usepackage[authoryear,round,longnamesfirst]{natbib}
\usepackage{fancyhdr}
\usepackage{hyperref}
\hypersetup{
  linktocpage=true,
  linkcolor  = MidnightBlue,
  citecolor  = MidnightBlue,
  urlcolor   = MidnightBlue,
  colorlinks = true,
}

\newcommand{\firmlist}{%
  \setlength{\itemsep}{0.5\itemsep}\setlength{\parskip}{0.5\parskip}}

\title{Folding a 3D Euclidean space\thanks{This is an expanded version of the article published in: R.~J.~Lang, M.~Bolitho and Z.~You (eds.), \textit{Origami$^7$ -- The proceedings from the 7th International Meeting on Origami in Science, Mathematics and Education}  (Tarquin), Vol.~2, pp.\ 331--346, 2018}}
\author{Jorge C. Lucero\thanks{Dept.\ Computer Science, University of Bras\'{i}lia, Brazil. E-mail: \href{mailto:lucero@unb.br}{lucero@unb.br} }} 
\date{\today} 
\begin{document}

\maketitle

\begin{abstract}
This paper considers an extension of origami geometry to the case of ``folding'' a three dimensional (3D) space along a plane. First, all possible incidence constraints between given points, lines and planes  are analyzed by using the geometry of reflections. Next, a set of 3D elementary fold operations is defined, which satisfy specific combinations of constraints with a finite number of solutions. The set consists of 47 valid fold operations, and solutions to some of them are explored to determine their number and conditions of existence.
\end{abstract}

\section{Introduction}
\label{sec:introduction}
Origami is the traditional Japanese art of creating figures by folding a sheet of paper. In modern days, the term has been extended to denote constructions techniques based on folding operations, and a number of applications to science and technology have been proposed, e.g.,  in aerospace and automotive technology \citep{Cipra2001}, civil engineering \citep{Filipov2015}, biology \citep{Mahadevan2005}, robotics \citep{Felton2014}, and other areas.  Its geometry has been studied through the so-called axioms of origami, which are elementary fold operations that satisfy specific incidence constraints between given points and lines on a plane. The axioms were introduced over three decades ago \citep{Justin1986} and, since then, have been expressed under a variety of forms \citep[e.g.,][]{Alperin2000,Alperin2006,Ghourabi2013,Lucero2017}. Analyses using the axioms usually regard the folded medium (i.e., the ``sheet'') as a 2D Euclidean plane, and extensions to the sphere and the hyperbolic plane have also been considered \citep{Kawasaki2011,Alperin2011}.   

What if, instead of folding a 2D medium, folds in higher dimensions are considered? This question has been considered from a topological point of view \citep{Robertson1978}; here, an axiomatic approach is adopted. In regular 2D origami, folding operations take place within a 3D space. If we increase one dimension, then we have a 3D medium folded within a 4D (or higher) hyperspace. Further, instead of folding along a line, the 3D folds take place along a plane. In this article, the 3D case is explored by extending the folding axioms to an Euclidean space and considering incidence constraints between points, lines and planes. Following a previous study \citep{Lucero2017}, the analysis is based on the geometry of reflections.

\section{Reflections in 3D}
\label{sec:section}
Throughout the analysis, points are denoted by capital letters ($P$, $Q$ etc.), lines by small Latin letters ($m$, $n$ etc.), and planes by small Greek letters ($\pi$, $\tau$ etc.) except the fold plane which is denoted by the special symbol $\Delta$. The notations $P\in m$ and $P\in \pi$ mean that point $P$ is on line $m$ and plane $\pi$, respectively. Further, $m\subset\pi$ means that line $m$ is contained in plane $\pi$. 

\subsection{Reflection of a point}
\begin{figure}
\centering
\includegraphics{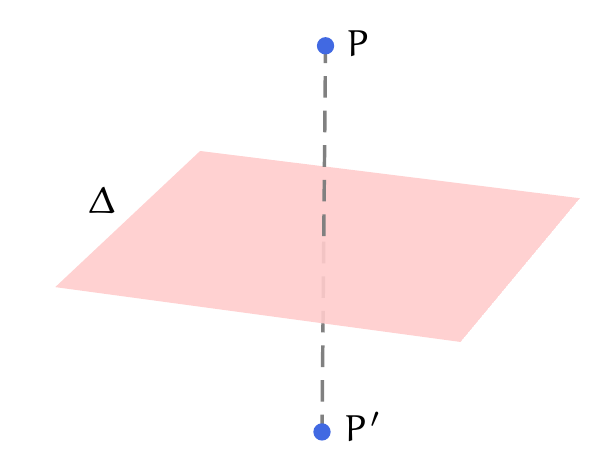}
\caption{Reflection of point $P$ in plane $\Delta$.} 
\label{refl1}
\end{figure}

Following the definition of reflection in a line \citep{Martin1998}, reflection in a plane is defined as  (see Fig.~\ref{refl1}):

\begin{definition}\label{defreflx} Given a plane $\Delta$, the reflection $\mathcal{F}_\Delta$ in $\Delta$ is the mapping on the set of points in the 3D space such that for point $P$
\[
\mathcal{F}_\Delta(P)=
\left\{
\begin{array}{ll}
P & \text{if $P\in \Delta$},\\
P' & \text{if $P\notin \Delta$ and $\Delta$ is the perpendicular bisector plane of segment $\overline{PP'}$}.
\end{array}
\right.
\]
\end{definition}

It is easy to see that $\mathcal{F}_\Delta(P)=P'$ iff $\mathcal{F}_\Delta(P')=P$.

\subsection{Reflection of a line}
Reflection of a line $m$ in $\Delta$ is obtained by reflecting every point in $m$. Therefore, $\mathcal{F}_\Delta(m)=\{\mathcal{F}_\Delta(P)|\,P\in m\}$. 

\begin{figure}
\centering
\includegraphics{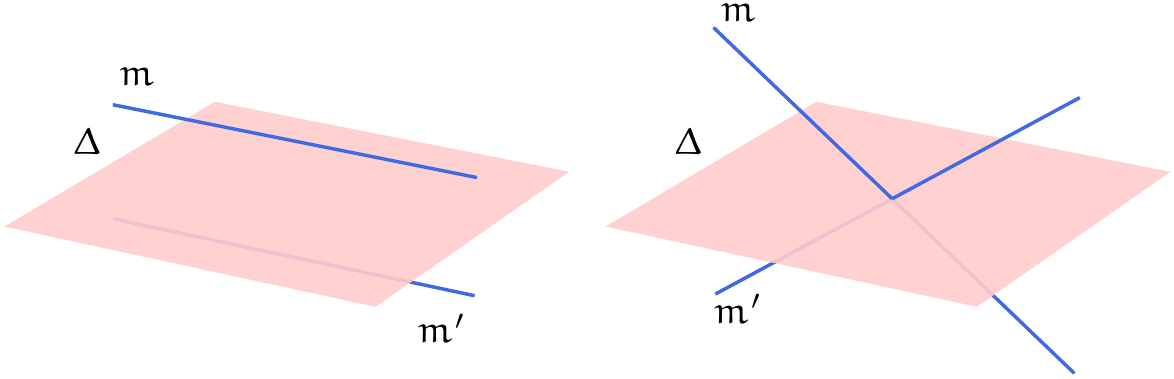}
\includegraphics{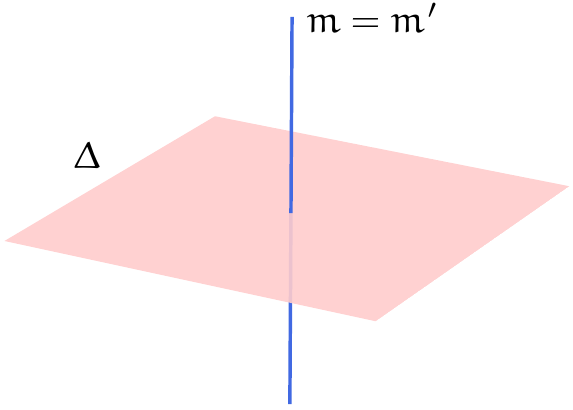}
\caption{Reflection of line $m$ in plane $\Delta$. Top left: $m\parallel\Delta$. Top right: $m\nparallel\Delta$. Bottom: $m\bot\Delta$.} 
\label{refl2}
\end{figure}

Let $m'=\mathcal{F}_\Delta(m)$. Then, lines $m$ and $m'$ are coplanar, and the plane defined by them is perpendicular to $\Delta$. 
Consider the following cases:
\begin{enumerate}
\firmlist
\item If $m$ and $\Delta$ are parallel ($m\parallel \Delta$), then $m\parallel m'$ (Fig.~\ref{refl2}, top left).
\item\label{it2} If $m$ and $\Delta$ are not parallel ($m\nparallel \Delta$)  then  the intersection of $\Delta$ with the plane spanned by $m$ and $m'$ is a bisector of the angle between $m$ and $m'$ (Fig.~\ref{refl2}, top right).
\item If $m\subset\Delta$, then every point $P\in m$ is its own reflection (i.e., $\mathcal{F}_\Delta(P)=P$), and therefore $m'=m$.
\item If $m$ and $\Delta$ are perpendicular ($m\bot\Delta$), then the reflection of every point of $m$ is also on $m$; therefore, $m=m'$ (Fig.~\ref{refl2}, bottom). Also, $\Delta$ divides $m$ into two half-lines (rays), and each half-line is reflected onto the other. Thus, for every point $P$ on $m$ not on the intersection with $\Delta$, $\mathcal{F}_\Delta(P)\neq P$. 
\end{enumerate}

\subsection{Reflection of a plane}
Reflection of a plane $\pi$ in $\Delta$ is similarly obtained by reflecting every point in $\pi$. Therefore, $\mathcal{F}_\Delta(\pi)=\{\mathcal{F}_\Delta(P)|\,P\in \pi\}$. 

\begin{figure}
\centering
\includegraphics{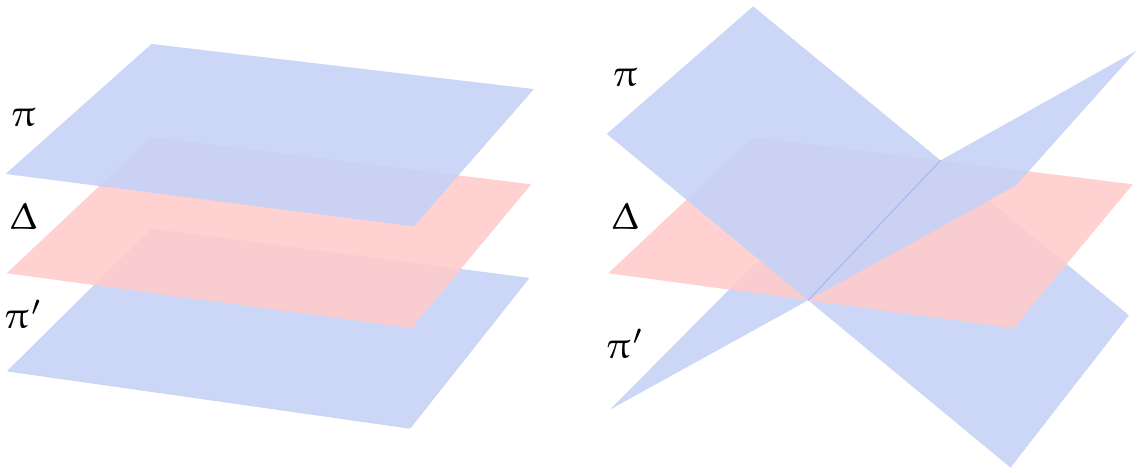}
\includegraphics{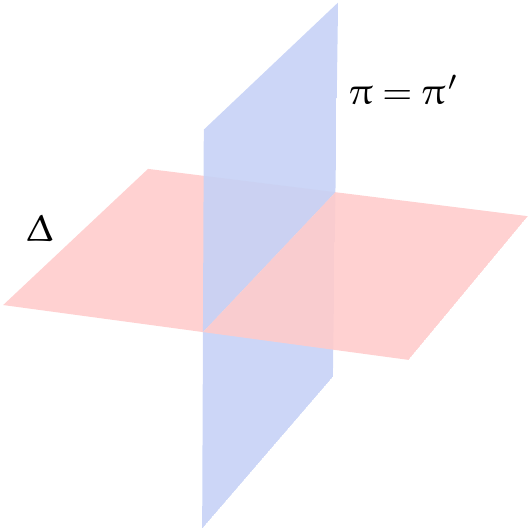}
\caption{Reflection of plane $\pi$ in plane $\Delta$. Top left: $\pi\parallel\Delta$. Top right: $\pi\nparallel\Delta$. Bottom: $\pi\bot\Delta$.} 
\label{refl3}
\end{figure}

Let $\pi'=\mathcal{F}_\Delta(\pi)$. Then:
\begin{enumerate}
\firmlist
\item If $\pi$ and $\Delta$ are parallel, then $\pi\parallel \pi'$ (Fig.~\ref{refl3}, top left).
\item If $\pi$ and $\Delta$ are not parallel  then  $\Delta$ is a bisector of the dihedral angle between $\pi$ and $\pi'$ (Fig.~\ref{refl3}, top right).
\item If $\pi=\Delta$, then every point $P\in\pi$ is its own reflection (i.e., $\mathcal{F}_\Delta(P)=P$), and therefore $\pi'=\pi$.
\item If $\pi$ and $\Delta$ are perpendicular, then the reflection of every point of $\pi$ is also on $\pi$; therefore, $\pi=\pi'$ (Fig.~\ref{refl3}, bottom). Also, $\Delta$ divides $\pi$ into two half-planes, and each half-plane is reflected onto the other. Thus, for every point $P$ on $\pi$ not on the intersection with $\Delta$, $\mathcal{F}_\Delta(P)\neq P$. 
\end{enumerate}

\section{Incidence constraints}
In 2D origami, elementary single-fold operations are defined in terms of incidence constraints between pairs of objects (points or lines) that must be satisfied with a fold \citep{Alperin2006, Ghourabi2013, Justin1986, Lucero2017}. Each constraint involves an object $\alpha$ and the image $\mathcal{F}_\Delta(\beta)$ of an object $\beta$ by reflection in the fold line $\Delta$. The symmetry of the reflection mapping implies that all incidence relations are symmetric.
 
Here, the same concept is extended to 3D: a plane is included in the objects that may be involved in a incidence constraint, and folding is performed along a plane instead of a line. 
A total of 12 different 3D incidences constraints may be defined, as follows (see also Table \ref{table0}). For convenience, incidences involving distinct objects (i.e., $\alpha\neq \beta$) are defined separately from those involving the same object (i.e., $\alpha= \beta$).

\subsection{Incidence \texorpdfstring{$I_1$: \normalfont$\mathcal{F}_\Delta(P)= Q$, with $P\neq Q$}{I1}}
\label{I1}
In this incidence, the reflection of a given point $P$ coincides with another given point $Q$. Its solution is the unique fold plane $\Delta$ which is the perpendicular bisector of segment $\overline{PQ}$ (Fig.~\ref{refl1}, with $Q=P'$.). 

The case $P=Q$ is treated in incidence $I_8$ (\S \ref{I8}).

\subsection{Incidence \texorpdfstring{$I_2$: \normalfont{$\mathcal{F}_\Delta(m)= n$, with $m\neq n$}}{I2}}
\label{I2}
\begin{figure}
\centering
\includegraphics{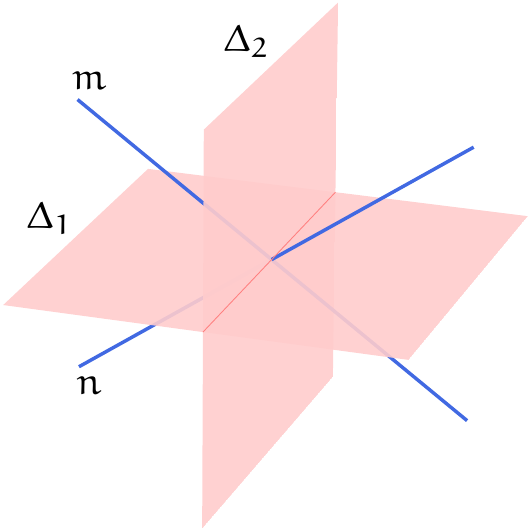}
\caption{Incidence $I_2$ in the case of $m\nparallel n$.}
\label{s2}
\end{figure}

In this incidence, the reflection of a given line $m$ coincides with another given line $n$. 
The following cases are possible:

\begin{enumerate}
\firmlist
\item If $m$ and $n$ are coplanar and $m\nparallel n$, there are two possible fold planes that satisfy the incidence. The fold planes are perpendicular between them and their intersections with the plane spanned by $m$ and $n$ are bisectors of the angles between those lines (Fig. \ref{s2}).  

\item If $m\parallel n$, there is only one solution, which is a fold plane perpendicular to the plane spanned by $m$ and $n$ and equidistant to them (Fig. \ref{refl2}, top left, with $n=m'$).
\item If $m$ and $n$ are not coplanar, then the incidence does not have a solution.
\end{enumerate}

The case $m=n$ is treated in incidences $I_9$ (\S \ref{I9}) and $I_{10}$ (\S \ref{I10}).

\subsection{Incidence \texorpdfstring{$I_3$: \normalfont{$\mathcal{F}_\Delta(m)\cap n\neq \emptyset$, with $m\cap n= \emptyset$}}{I3}}
\label{I3}
In this incidence, the reflection of a given line $m$ intersects another line $n$, and the case in which $m$ and $n$ already intersect is excluded. 

\begin{figure}
\centering
\includegraphics{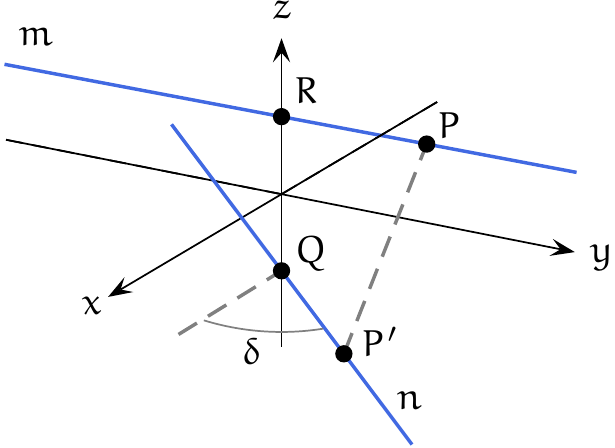}
\caption{Diagram for incidence $I_3$ in the case of $m$ and $n$ not coplanar.}
\label{i2b}
\end{figure}

Assume first that $m$ and $n$ are not coplanar and consider a Cartesian system of coordinates $x, y, z$ such that $m$ lies in the $yz$ coordinate plane and passes through point $R(0, 0, 1)$, and $n$ lies in the  $z=-1$ plane and passes through point $Q(0, 0, -1)$. (Fig.~\ref{i2b}). Let $P(0, t, 1)$ be an arbitrary point in $m$, where $t$ is a parameter, and its reflection in $n$ be $P'(s\cos\delta, s\sin\delta, -1)$, where $\delta\neq \pi/2$ is the angle with the $x$-axis, and $s$ is a parameter. 

 A fold plane $\Delta$ may be described by a vector equation of the form $\mathbf{n}\cdot\mathbf{x}=\mathbf{n}\cdot\mathbf{x}_0$, where $\mathbf{n}$ is a normal vector, $\mathbf{x}=(x, y, z)$, and $\mathbf{x}_0$ is a point on $\Delta$. A normal vector for $\Delta$ is $\overrightarrow{PP'}=(s\cos\delta, s\sin\delta - t , -2)$, and $\Delta$ passes through the midpoint of $\overline{PP'}$ at $(s\cos\delta, t+s\sin\delta, 0)/2$. Then, $\Delta$ has an equation
\begin{equation}
sx\cos\delta+y(s\sin\delta-t) -2z = \frac{s^2-t^2}{2}.
\label{deltai2b}
\end{equation}

Eq.~(\ref{deltai2b}) defines a family of fold planes in two parameters ($t$ and $s$).  The envelope of the family may be found as explained in Appendix \ref{ApA}, which results in
\begin{equation}
x^2\cos^2\delta+2xy\sin\delta\cos\delta -y^2\cos^2\delta-4z=0.
\label{envelopeb}
\end{equation}

Eq.~(\ref{envelopeb}) represents a quadric surface, and it may be put into normal form by rotating the coordinate system around the $z$-axis so as to eliminate the term in $xy$. 
A rotation by an angle $\delta/2$
produces
\begin{equation}
u^2\cos\delta-v^2\sin\delta -4z=0,
\label{envelopeb1}
\end{equation}
where
\begin{equation}
\begin{pmatrix}
u\\
v
\end{pmatrix}
=
\begin{pmatrix}
\cos\delta/2 & \sin\delta/2\\
-\sin\delta/2 &\cos\delta/2
\end{pmatrix}
\begin{pmatrix}
x\\
y
\end{pmatrix}.
\end{equation}

\begin{figure}
\centering
\includegraphics{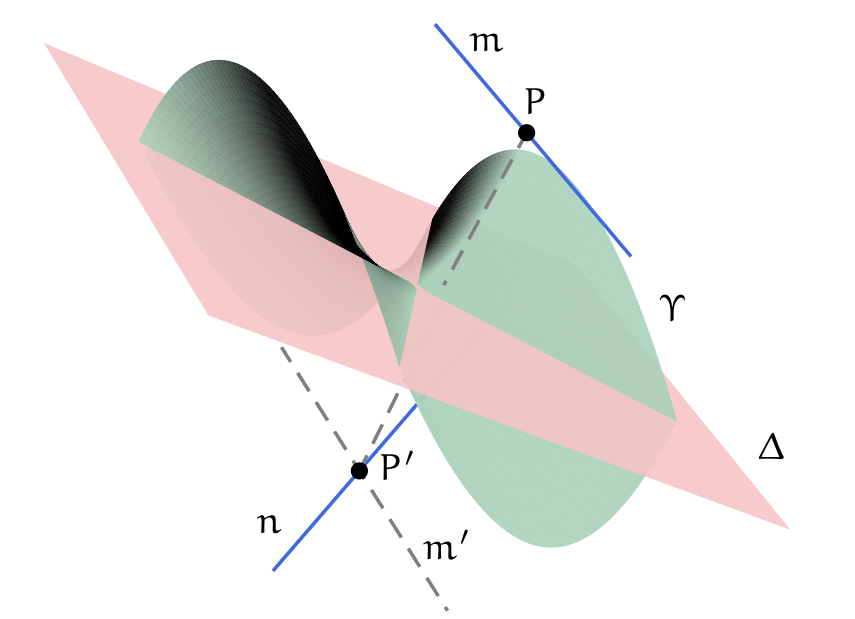}
\caption{Example of fold plane for incidence $I_3$ in the case that $m$ are not coplanar. Plane $\Delta$ reflects $m$ onto $m'$, which intersects $n$ at $P'$, and is tangent to the hyperbolic paraboloid $\Upsilon$.}
\label{hyper}
\end{figure}

Eq.~(\ref{envelopeb1}) describes an hyperbolic paraboloid that opens up along the $u$ axis and down along the $v$ axis. An example is shown in Fig.~\ref{hyper}. 

If $m$ and $n$ are coplanar and do not intersect, then they are parallel. Letting  $\delta=\pi/2$ in  Eq.~(\ref{deltai2b}) produces
\begin{equation}
y(t-s) +2z = \frac{t^2-s^2}{2}.
\label{deltai2c}
\end{equation}
This equation describes any plane parallel to the $x$-axis that intersects the $y$-axis (including the case that the fold plane is the $xy$ plane itself) and is not perpendicular to it. In fact, note that $z=0$ implies $t=s$ and therefore the fold plane is $z=0$ (i.e., the $xy$ plane), or $y=(t+s)/2$ and therefore the fold plane intersect the $y$ axes. Also, the coefficient for the term in $z$ is not zero, so the fold plane cannot be parallel to the $z$-axis. 

If $m$ and $n$ already intersect at a point $P$ and $m\neq n$, then there are two families of possible fold planes:
\begin{enumerate} 
\firmlist
\item Any plane $\Delta$ that passes through $P$. Such a plane reflects $P$ onto itself, and therefore $m'$and $n$ intersect at the same point $P$. This case may be included within incidence $I_7$ (\S \ref{I7}).
\item Any plane $\Delta$ that is perpendicular to the plane spanned by $m$ and $n$ and either $\Delta$ contains the angle bisector of $m$ and $n$, denoted as $\ell$, or $\Delta$  does not contain $\ell$ and is not parallel to it. Note that, if $\Delta$ does not contain $\ell$ but is parallel to it, then $m'$ results parallel to $n$ and therefore both lines do not intersect. In any other case, $m'$ and $n$ are coincident (if $\Delta$ contains $\ell$) or they intersect at some point. 

This case may be included within incidence $I_{10}$ (\S \ref{I10}), which is solved by any fold plane perpendicular to a given plane $\pi$. The plane $\pi$ is the plane spanned by $m$ and $n$, and the case in which the fold plane is parallel to $\ell$ must be excluded.  
\end{enumerate} 

If $m=n$ then any fold plane solves the incidence, with the exception of planes parallel to $m$ and not containing it. This case  does not constitute a valid constraint (the family of fold planes still has three degrees of freedom; see \S \ref{efo}) and may be disregarded.

\subsection{Incidence \texorpdfstring{$I_4$: \normalfont{$\mathcal{F}_\Delta(\pi)= \tau$, with $\pi\neq \tau$}}{I4}}
\label{I4}
\begin{figure}[!htb]
\centering
\includegraphics{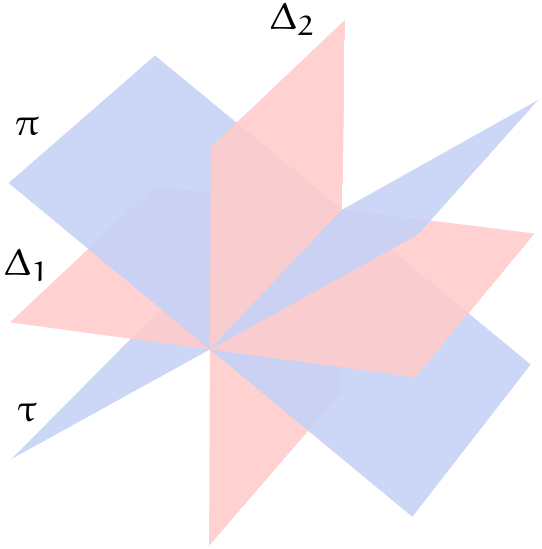}
\caption{Incidence $I_4$ in the case of $\pi\nparallel \tau$. } 
\label{s3}
\end{figure}

In this incidence, the reflection of a given plane $\pi$ coincides with another given plane $\tau$. Two cases are possible:

\begin{enumerate}
\firmlist
\item If  $\pi\nparallel\tau$, there are two possible fold planes that satisfy the incidence. The fold planes are perpendicular and they  bisect the dihedral angles between $\pi$ and $\tau$ (Fig. \ref{s3}).  
\item If $\pi\parallel \tau$, there is only one solution, which is a fold plane parallel and equidistant to both $\pi$ and $\tau$ (Fig. \ref{refl3}, top left, letting $\tau=\pi'$).
\end{enumerate}

The case $\pi=\tau$ is treated in incidences $I_{11}$ (\S \ref{I11}) and $I_{12}$ (\S \ref{I12}).

\subsection{Incidence \texorpdfstring{$I_5$: \normalfont$\mathcal{F}_\Delta(P)\in m$, with $P\notin m$}{I5}}
\label{I5}
In this incidence, the reflection of a given point $P$ is on a given line $m$, and the case in which $P$ is already on $m$ is excluded. In 2D folding, reflecting a point $P$ on a line $m$ produces fold lines that are tangents to a parabola with focus $P$ and directrix $m$ \citep{Alperin2000, Martin1998}. The 3D case is similar, but a parabolic cylinder is produced instead.

\begin{figure}
\centering
\includegraphics{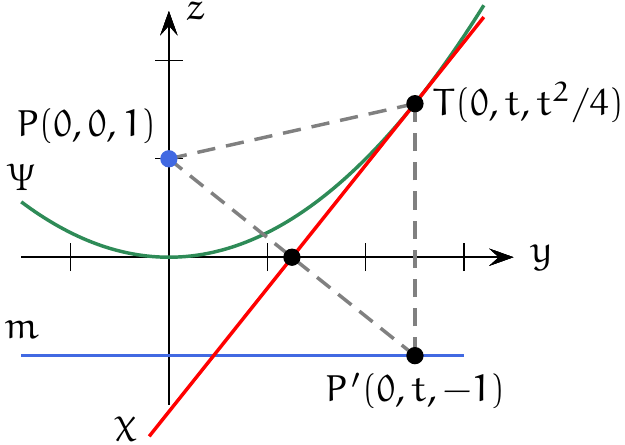}
\caption{Diagram for incidence $I_5$. Line $\chi$ is the intersection of the fold plane $\Delta$ with the plane $yz$.} 
\label{parab}
\end{figure}

Choose a Cartesian system of coordinates $x, y, z$ so that point $P$ is at $(0, 0, 1)$, and line $m$ passes through point $(0, 0, -1)$ and is parallel to the $y$-axes. Let $P'=\mathcal{F}_\Delta(P)\in m$ be at $(0, t, -1)$, where $t$ is a free parameter (Fig.~\ref{parab}).  

A normal vector for the fold plane $\Delta$ is $\overrightarrow{PP'}=(0, t, -2)$, and $\Delta$ passes through the midpoint of $\overline{PP'}$ at $(0, t/2, 0)$. Thus, an equation for $\Delta$ is
\begin{equation}
2ty-4z = t^2,
\label{deltaplane}
\end{equation}
which defines a one-parameter family of planes perpendicular to the $yz$ coordinate plane. 
The envelope of the family is given by
\begin{equation}
z=\frac{y^2}{4}.
\label{pp}
\end{equation}

\begin{figure}
\centering
\includegraphics{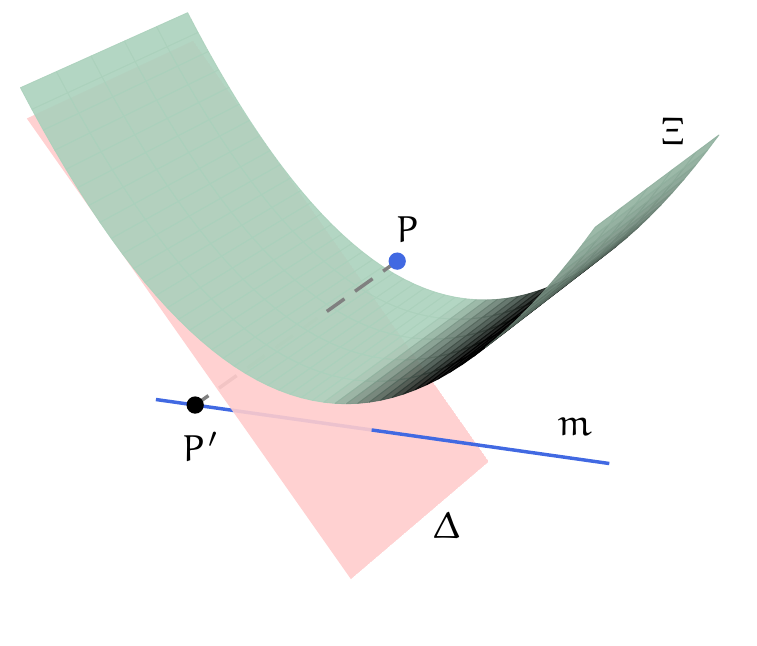}
\caption{Example of fold plane for incidence $I_5$. Plane $\Delta$ reflects $P$ onto $P'\in m$ and is tangent to the parabolic cylinder $\Xi$.}
\label{cylp}
\end{figure}

Eq.~(\ref{pp}) describes a parabolic cylinder generated by a parabola in the $yz$ plane, with focus at $P$ and directrix $m$, when translated in direction parallel to the $x$ axis.  An example is shown in Fig.~\ref{cylp}.

If $P\in m$, then any plane containing $P$ or perpendicular to $m$ reflects $P$ onto $m$. Those two cases are contemplated by incidences $I_8$ (\S \ref{I8}) and $I_9$ (\S \ref{I9}), respectively.

\subsection{Incidence \texorpdfstring{$I_6$: \normalfont{$\mathcal{F}_\Delta(P)\in \pi$, with $P\notin \pi$}}{I6}}
\label{I6}
In this incidence, the reflection of a given point $P$ is on a given plane $\pi$, and the case in which $P$ is already on $\pi$ is excluded.

The analysis is similar to the previous incidence. In the same Cartesian system, let plane $\pi$ be described by $z=-1$, and $P'=\mathcal{F}_\Delta(P)\in \pi$ be located at $(s, t, -1)$, where $s$ and $t$ are free parameters. Then, the same diagram in Fig.~\ref{parab} applies, where $m$ is now the intersection line of $\pi$ with the $yz$ plane. 

A normal vector for the fold plane $\Delta$ is $\overrightarrow{PP'}=(s, t, -2)$, and $\Delta$ passes through the midpoint of $\overline{PP'}$ at $(s/2, t/2, 0)$. Thus, an equation for $\Delta$ is 
\begin{equation}
2sx+2ty-4z = s^2+t^2.
\label{deltaplane1}
\end{equation}
This equation represents a family of planes with two parameters ($s$ and $t$), and the family has an envelope given by
\begin{equation}
x^2+y^2-4z=0.
\label{parabld}
\end{equation}

\begin{figure}[!htb]
\centering
\includegraphics{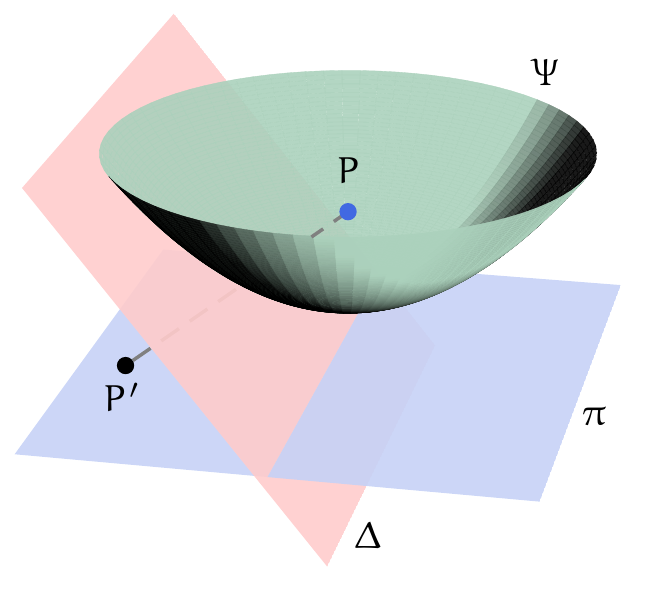}
\caption{Example of fold plane for incidence $I_6$. Plane $\Delta$ reflects $P$ onto $P'\in \pi$ and is tangent to the paraboloid $\Psi$.}
\label{paraboloid}
\end{figure}

Eq. (\ref{parabld}) describes a paraboloid generated by rotation around the $z$ axis of a parabola with vertex at $(0, 0, 0)$ and axis coincident with the $z$ axis (Fig.~\ref{paraboloid}). The parabola has focus at $P$ and its directrix is any line in $\pi$ that intersects the $z$ axis.  An example is shown in Fig.~\ref{paraboloid}.

If $P\in \pi$, then any plane containing $P$ or perpendicular to $\pi$ reflects $P$ onto $\pi$. Those two cases are contemplated by incidences $I_8$ (\S \ref{I8}) and $I_{11}$ (\S \ref{I11}), respectively.

\subsection{Incidence \texorpdfstring{$I_7$: \normalfont$\mathcal{F}_\Delta(m)\subset \pi$, with $m\not\subset \pi$}{I7}}
\label{I7}
In this incidence, the reflection of a given line $m$ is on a given plane $\pi$, and the case in which $m$ is already on $\pi$ is excluded.

\begin{figure}
\centering
\includegraphics{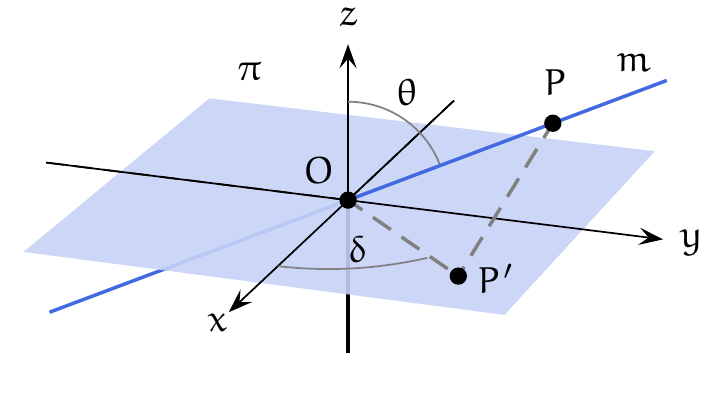}
\caption{Diagram for incidence $I_7$ in the case of $m\nparallel\pi$.}
\label{i61}
\end{figure}

Assume first $m\nparallel\pi$, and choose a Cartesian system of coordinates $x, y, z$ so that the point of intersection of $m$ and $\pi$, denoted by $O$, is located at $(0, 0, 0)$, plane $\pi$ is the $xy$ coordinate plane, and $m$ is a line in $yz$ plane forming an angle $\theta\neq \pi/2$ with the $z$-axis (Fig.~\ref{i61}). An arbitrary point $P\in m$ is located at $(0, \sin\theta, \cos\theta)t$, where $t$ is a parameter. Its reflection $P'=\mathcal{F}_\Delta(P)$ 
is in $\pi$ at $(\cos\delta, \sin\delta, 0)t$, where $\delta$ is the angle with the $x$-axis.

A normal vector for the fold plane $\Delta$ is $\overrightarrow{PP'}=(\cos\delta, \sin\delta-\sin\theta, -\cos\theta)t$, and $\Delta$ passes through $O$. Thus, an equation for $\Delta$ is 
\begin{equation}
x\cos\delta+y(\sin\delta-\sin\theta)-z\cos\theta=0.
\label{deltaplane2}
\end{equation}
This equation defines a family of fold planes in one parameter ($\delta$), which has an envelope given by
\begin{equation}
x^2+y^2\cos^2\theta-2yz\sin\theta\cos\theta-z^2\cos^2\theta=0.
\label{envelope3}
\end{equation}

Eq.~(\ref{envelope3}) describes a quadric surface, and the equation may be put into normal form by rotating the coordinate system around the $x$-axis so as to eliminate the term in $yz$. A rotation by an angle $-\theta/2$ 
produces
\begin{equation}
\frac{x^2}{\cos\theta}+u^2-v^2=0,
\label{envelope4}
\end{equation}
where
\begin{equation}
\begin{pmatrix}
u\\
v
\end{pmatrix}
=
\begin{pmatrix}
\cos\theta/2 & -\sin\theta/2\\
\sin\theta/2 &\cos\theta/2
\end{pmatrix}
\begin{pmatrix}
y\\
z
\end{pmatrix}.
\end{equation}

\begin{figure}
\centering
\includegraphics{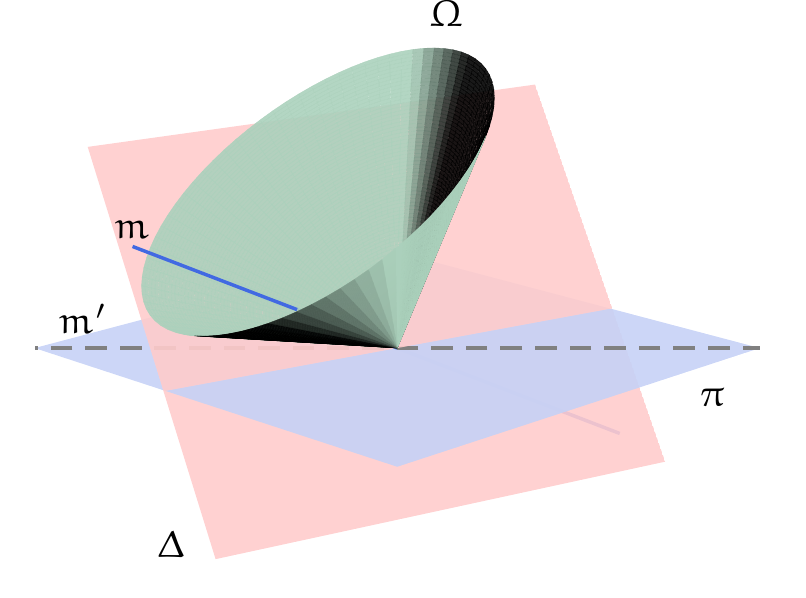}
\caption{Example of fold plane for incidence $I_7$ in the case of $m\nparallel\pi$. Plane $\Delta$ reflects $m$ onto $m'\in \pi$ and is tangent to the elliptical cone $\Omega$. For clarity of the figure, only the upper nappe of $\Omega$ is shown.}
\label{cone2}
\end{figure}

Eq.~(\ref{envelope4}) describes an elliptical cone that opens along the $v$-axis. An example is shown in Fig.~\ref{cone2}.

\begin{figure}
\centering
\includegraphics{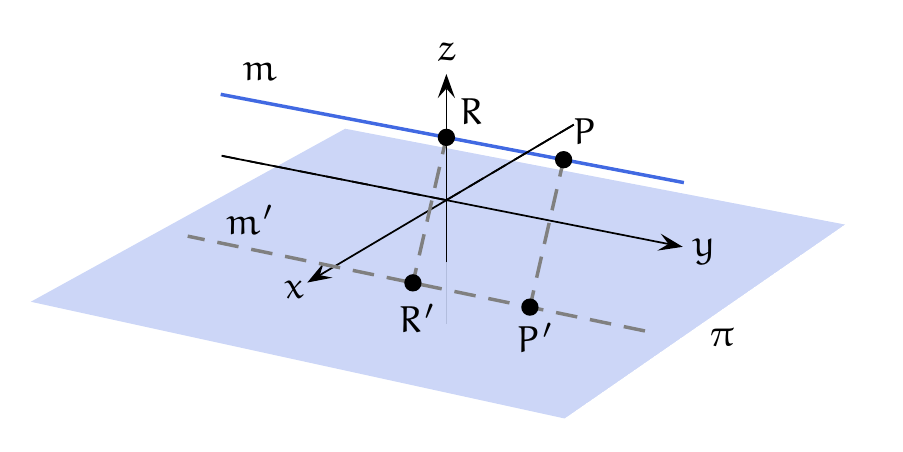}
\caption{Diagram for incidence $I_7$ in the case of $m\parallel\pi$.}
\label{i62}
\end{figure}

Next, consider the case $m\parallel\pi$. 
Let $\pi$ be the plane $z=-1$, and $m$ be a line in the $yz$ coordinate plane passing through point $R$ at $(0, 0, 1)$ and parallel to the $y$-axis (Fig.~\ref{i62}). The reflection of $m$ is $m'=\mathcal{F}_\Delta(m)\subset \pi$, with $m'\parallel m$,  and passes through $R'=\mathcal{F}_\Delta(R)$ at $(k, 0, -1)$, where $k$ is a parameter.

An arbitrary point $P\in m$ is located at  $(0, t, 1)$, where $t$ is a parameter, and its reflection $P'=\mathcal{F}(P)$  is at $(k, t, -1)$. Following similar steps as before, we find a fold plane $\Delta$ with equation
\begin{equation}
-2kx+4z=-k^2.
\label{deltaplane3}
\end{equation}

\begin{figure}
\centering
\includegraphics{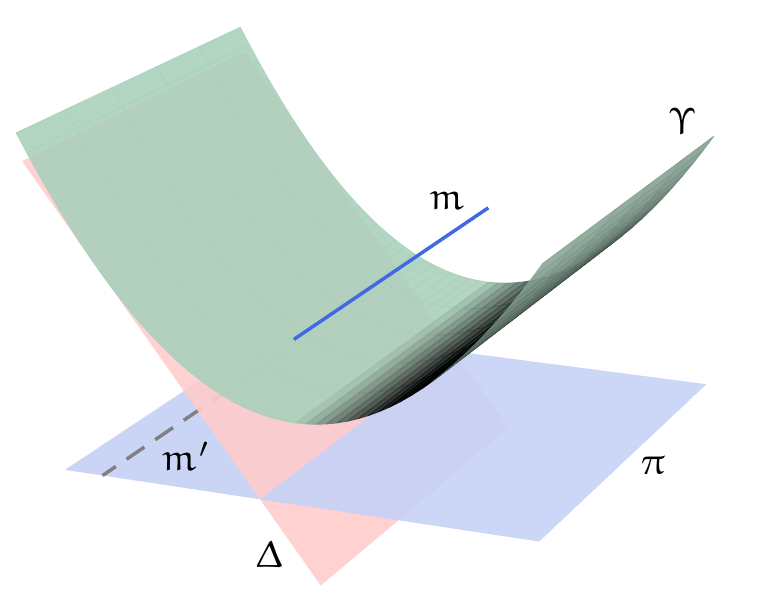}
\caption{Example of fold plane for incidence $I_7$ in the case of $m\parallel\pi$. Plane $\Delta$ reflects $m$ onto $m'\in \pi$ and is tangent to the parabolic cylinder $\Upsilon$.}
\label{cone3}
\end{figure}

Eq.~(\ref{deltaplane3}) defines a family of fold planes in parameter $k$, and its envelope is given by
\begin{equation}
x^2-4z=0,
\label{envelope5}
\end{equation}
which describes a parabolic cylinder generated by a parabola in the $xz$ plane when translated in direction parallel to the $y$-axis. The parabola has focus in $R$ and its directrix is a line parallel to the $x$-axis and passing through the point $(0, 0, -1)$.  An example is shown in Fig.~\ref{cone3}.

If $m\subset \pi$, then any plane containing $m$ or perpendicular to $\pi$ reflects $m$ onto $\pi$. Those two cases are contemplated by incidences I10 (\S \ref{I10}) and I11 (\S \ref{I11}), respectively.

\subsection{Incidence \texorpdfstring{$I_8$: \normalfont$\mathcal{F}_\Delta(P)=P$}{I8}}
\label{I8}
In this incidence, the reflection of a given point $P$ coincides with itself. It is satisfied by any fold plane $\Delta$ passing through $P$. An arbitrary normal direction for $\Delta$ may be defined, e.g., by the azimuthal angle $\theta$ in the $xy$ and the polar angle $\phi$ from the $z$-axis. Therefore, the solution to the incidence is a family of fold planes with two free parameters, namely, $\theta$ and $\phi$.

\subsection{Incidence \texorpdfstring{$I_9$: \normalfont$\mathcal{F}_\Delta(m)=m$, and $\exists P\in m, \mathcal{F}_\Delta(P)\neq P$}{I9}}
\label{I9}
Both this and the next incidence consider the reflection of a line $m$ to itself. In the current incidence, half of $m$, defined from an arbitrary point $R\in m$, is reflected upon the opposite half. 

The position of point $R$ may be specified by its distance $s$ from a particular point $P_0\in m$. A fold plane is perpendicular to $m$ and passes through point $R$. Therefore, the solution to the incidence is a family of fold planes with one free parameter ($s$).

\subsection{Incidence \texorpdfstring{$I_{10}$: \normalfont$\mathcal{F}_\Delta(m)=m$, and $\forall P\in m, \mathcal{F}_\Delta(P)= P$}{I10}}
\label{I10}
This is the second case of reflection of line $m$ to itself. In this incidence, each point $P\in m$ is reflected to itself, and therefore the incidence is satisfied by a fold plane containing $m$. There is an infinite number of such planes, an all of them have $m$ as intersection. Each plane may be  specified by the dihedral angle $\varphi$ from a reference plane; therefore, the solution to the incidence is a family of fold planes with one free parameter ($\varphi$).

\subsection{Incidence \texorpdfstring{$I_{11}$: \normalfont$\mathcal{F}_\Delta(\pi)=\pi$, and $\exists P\in \pi, \mathcal{F}_\Delta(P)\neq P$}{I11}}
\label{I11}
Both this and the next incidence consider the reflection of a plane $\pi$ to itself. In the current case, half of $\pi$, defined from an arbitrary line $m\subset \pi$, is reflected upon the opposite half. Any fold plane perpendicular to $\pi$ satisfies the incidence, and it may be specified by an equation of the form $\mathbf{n}\cdot\mathbf{x}=k$, where $\mathbf{n}$ is a normal vector and $k$ a parameter. 
Vector $\mathbf{n}$ is parallel to $\pi$, and an arbitrary direction may be specified by its angle $\delta$ from a reference  direction parallel to $\pi$. 
Therefore, the solution to the incidence is a family of fold planes with two free parameter, namely, $k$ and $\delta$.

\subsection{Incidence \texorpdfstring{$I_{12}$: \normalfont$\mathcal{F}_\Delta(\pi)=\pi$, and $\forall P\in \pi, \mathcal{F}_\Delta(P)= P$}{I12}}
\label{I12}
This is the second case of reflection of plane $\pi$ to itself. In this case, each point $P\in \pi$ is reflected to itself, and therefore the fold plane is $\Delta=\pi$

\section{Elementary fold operations}
\subsection{Definition}
\label{efo}
 
A plane in 3D space is an object with three degrees of freedom.\footnote{A plane $\Delta$ may be defined by an equation of the form $ax+by+cz+d=0$, where $a$, $b$, $c$ and $d$ are constants, and $(a, b, c)$ is a normal vector to $\Delta$. A vector in arbitrary direction may be defined by  a polar angle $0\le \theta\le \pi$ and an azimuthal angle $0\le \varphi< 0$, and letting $a=\sin\theta\cos \varphi$, $b=\sin\theta\sin\varphi$, and $c=\cos\theta$.  Therefore, three parameters must be set in order to define any fold plane, namely, $d$, $\theta$ and $\varphi$.} 
When an incidence constraint is set for the fold plane $\Delta$, satisfying the constraint consumes a number of degrees of freedom, and that number is called the codimension of the constraint. Incidences $I_1$, $I_2$, $I_4$ and $I_{12}$ have either a unique or a finite number of solutions; therefore, they define constraints of codimension 3. Incidences  $I_5$, $I_7$, $I_9$ and $I_{10}$ have a family of solutions with one free parameter and therefore they define constraints of codimension 2.  Finally, incidences  $I_3$, $I_6$, $I_8$ and $I_{11}$ have a family of solutions with two free parameter and therefore they define constraints of codimension 1.  Table \ref{table0} lists all the incidence constraints with their respective codimensions.

\begin{table}
\centering
\caption{Incidence constraints. $P$ and $Q$ are points; $m$ and $n$ are lines; $\pi$ and $\tau$ are planes.}
  \renewcommand{\arraystretch}{1.2}
\begin{tabular}{clc}
\hline
Incidence & Definition & Codimension\\
\hline
I1  & $\mathcal{F}_\Delta(P)=Q$, with $P\neq Q$ & 3\\
I2  & $\mathcal{F}_\Delta(m)=n$, with $m\neq n$ & 3 \\
I3  & $\mathcal{F}_\Delta(m)\cap n\neq \emptyset$, with $m\cap n= \emptyset$ & 1\\
I4  & $\mathcal{F}_\Delta(\pi)=\tau$, with $\pi\neq \tau$ & 3 \\
I5  & $\mathcal{F}_\Delta(P)\in m$, with $P\notin m$ & 2 \\
I6  & $\mathcal{F}_\Delta(P)\in \pi$, with $P\notin \pi$ & 1 \\
I7  & $\mathcal{F}_\Delta(m)\subset \pi$, with $m\not\subset \pi$ & 2 \\
I8  & $\mathcal{F}_\Delta(P)=P$ & 1\\
I9  & $\mathcal{F}_\Delta(m)=m$, and $\exists P\in m$, $\mathcal{F}_\Delta(P)\neq P$ &2\\
I10 & $\mathcal{F}_\Delta(m)=m$, and $\forall P\in m, \mathcal{F}_\Delta(P)= P$ & 2\\
I11 & $\mathcal{F}_\Delta(\pi)=\pi$, and $\exists P\in \pi$, $\mathcal{F}_\Delta(P)\neq P$ &1\\
I12 & $\mathcal{F}_\Delta(\pi)=\pi$, and $\forall P\in \pi, \mathcal{F}_\Delta(P)= P$ & 3\\
\hline
\end{tabular}
	\label{table0}
\end{table}

An elementary fold operation may be defined as the resolution of a minimal set of incidence constraints with a finite number of solutions \citep{Alperin2006}.  Let us recall that the total codimension of the combination is \textit{at most} the sum of the individual codimensions. Therefore, that sum must not be smaller than 3. 

Each of the incidences $I_1$, $I_2$, $I_4$ and $I_{12}$ already define an elementary operation. Here, it may be argued that incidence $I_{12}$ should be disregarded because it does not create a new plane. However, completeness of the set of possible incidences and operations demands its inclusion \citep[see a similar discussion for the case of 2D folding by][]{Lucero2017}. The other constraints must be applied in combinations, as follows:
\begin{enumerate}
\item Two constraints of respective codimensions 1 and 2. There are four constraints of each type, and so the total number of combinations is $4\times4=16$.
\item Two constraints both of codimension 2. The number $N$ of possible combinations when selecting $r$ objects from a set of $n$ objects and allowing repetitions is given by\cite{Rosen2012} 
\begin{equation}
N=\binom{n+r-1}{r}= \frac{(n+r-1)!}{r!(n-1)!}.
\label{binom}
\end{equation}
In the present case, two constraints ($r=2$) are selected from a total of four ($n=4$), which results in 10 combinations.
\item  Three constraints of codimension 1. Again, using Eq. ({\ref{binom}}) with $n=4$ and $r=3$ produces a number of 20 combinations.
\end{enumerate}

However, incidence $I_{11}$ may not be applied three times. If it is, then the fold plane $\Delta$ must be perpendicular to three given planes. If two of them are parallel, then the associated constraints are mutually redundant and therefore the total codimension is less than 3, with an infinite number of solutions. If all given planes are mutually nonparallel, then a fold plane does not exist (in an Euclidean space). Hence, this combination does not define a valid folding operation. 

For similar reasons,  incidence $I_{11}$ may not be combined with incidence $I_9$, and also incidence $I_9$ may not be used twice. 

The remaining combinations may not necessarily result in distinct operations. Thus, a total of at most 47 valid operations may be defined, and some of them are discussed in the next subsections. A notation in the format $a_1I_{b_1}+\cdots+a_kI_{b_k}$ is used, meaning that incidences $I_{b_1},\ldots, I_{b_k}$ are combined, and each one is used $a_1,\ldots, a_k$ times, respectively.

\subsection{Operation \texorpdfstring{$I_1+I_6$}{I5 + I6}}
\label{O6}
This operation combines one constraint of codimension 2 ($I_5$) and one of codimension 1 ($I_6$), and it may be stated as: \textit{Given a line $m$, a plane $\pi$, a point $P$ not on $m$, and a point $Q$ not on $\pi$, fold along a plane to place $P$ onto $m$ and $Q$ onto $\pi$}. Its solution is a fold plane tangent to the parabolic cylinder generated by $P$ and $m$ (Fig.\ \ref{cylp}) and the paraboloid generated by $Q$ and $\pi$ (Fig.\ \ref{paraboloid}, with $Q$ in place of $P$).  

Assume $P$ and $m$ as in \S \ref{I5}, point $Q$ at $(x_q, y_q, z_q)$ and $Q'$ at $(x'_q, y'_q, z'_q)$. The fold plane $\Delta$ is given by Eq. (\ref{deltaplane}) and it must be normal to segment $\overline{QQ'}$. Therefore, the normal vector of $\Delta$ is parallel to $\overline{QQ'}$ and satisfies
\begin{equation}
(0, t, -2) = k(x'_q-x_q, y'_q-y_q, z'_q-z_q),
\label{o61}
\end{equation}
where $k$ is a constant. Further, $\Delta$ passes through the midpoint of $\overline{QQ'}$ located at $((x_q +x'_q)/2, (y_q +y'_q)/2, (z_q+z'_q)/2)$. Replacing these coordinates into Eq. (\ref{deltaplane}) produces
\begin{equation}
t(y_q+y'_q)-2(z_q+z'_q)=t^2.
\label{o62}
\end{equation}

Eliminating $t$ and $k$ from Eqs. (\ref{o61}) and (\ref{o62}) gives
\begin{equation}
x_q=x'_q,
\label{o63a}
\end{equation}
and
\begin{equation}
2(y_q-y'_q)^2=-(y_q^2-{y'}_q^2)(z_q-z'_q)-(z_q^2-{z'}_q^2)(z_q-z'_q).
\label{o63}
\end{equation}

For a given plane $\pi$ described by $ax+by+cz+d=0$, where $a$, $b$, $c$ and $d$ are constants, the coordinates of $Q'$ satisfy
\begin{equation}
ax_q+by'_q+cz'_q+d=0
\label{o64}
\end{equation}
where Eq.~(\ref{o63a}) has already been used. 

\begin{figure}
\centering
\includegraphics{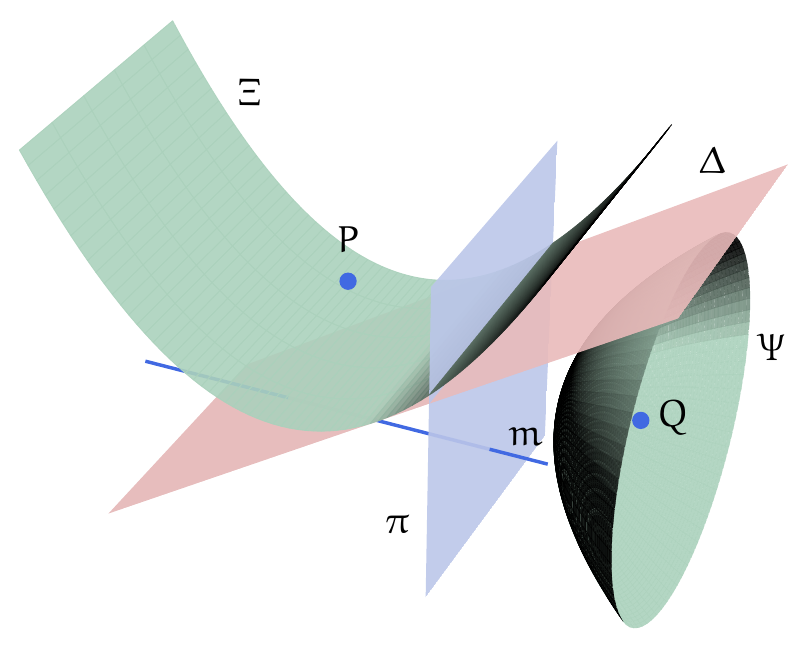}\\[-2ex]
\includegraphics{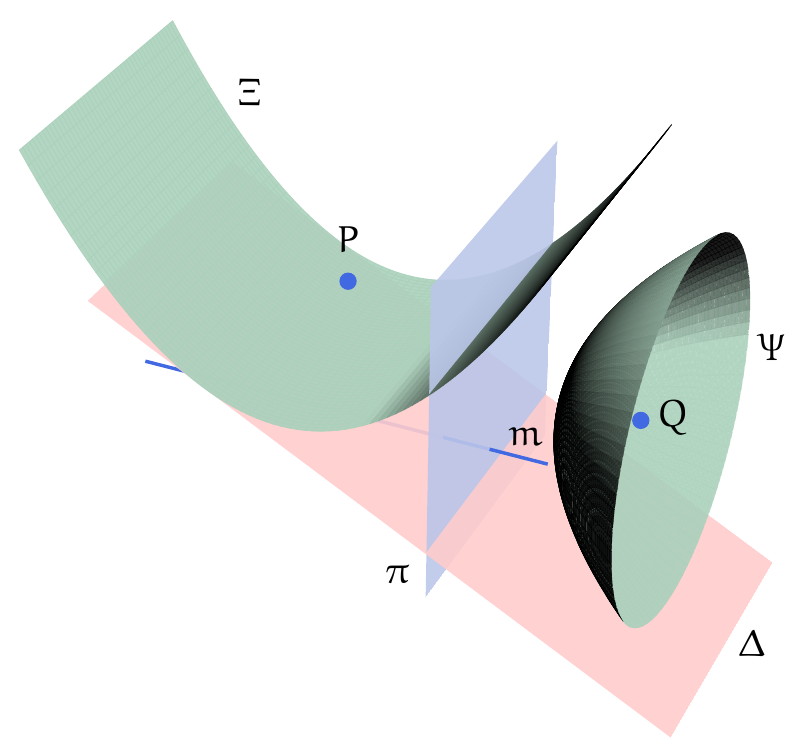}\\[-2ex]
\includegraphics{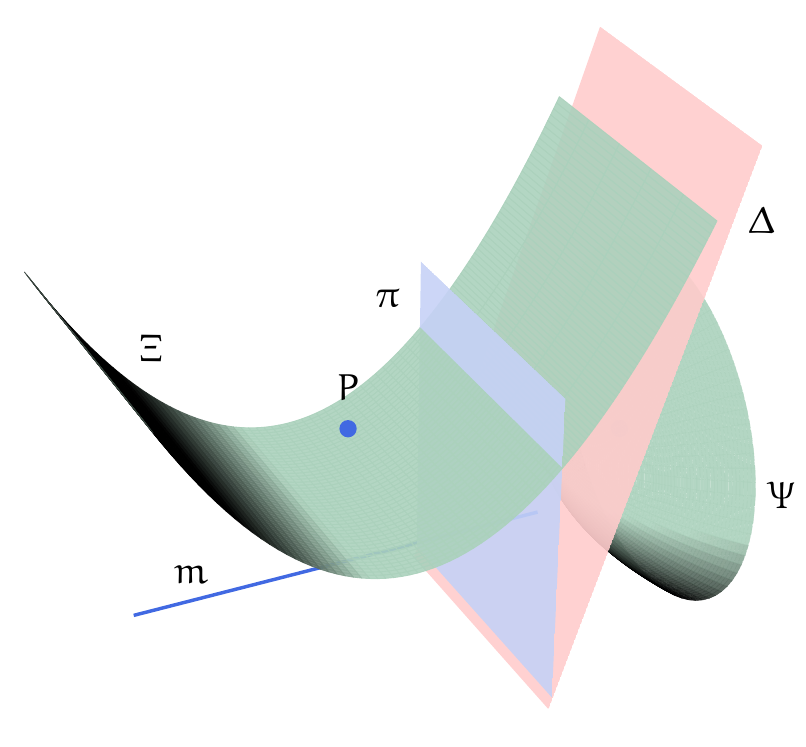}
\caption{Example of operation $I_5 + I_6$ with three solutions (the bottom figure has been rotated to improve clarity). The fold plane $\Delta$ is tangent to the parabolic cylinder $\Xi$, with focus $P$ and directrix line $m$, and the paraboloid $\Psi$, with focus $Q$ and directrix plane $\pi$.}
\label{I5+I6}
\end{figure}

Eqs.~(\ref{o63}) and (\ref{o64}) may be solved for $y'_q$ and $z'_q$. Three cases are possible:
\begin{enumerate}
\firmlist
\item If $P=Q$ then $y_q=0$, $z_q=1$. Further, if $\pi$ is the plane $z=-1$ (i.e., $\pi$ contains $m$ and is perpendicular to the plane spanned by $P$ and $m$), then $z'_q=-1$. In this case,  Eq.~(\ref{o63}) is satisfied by any value of $y'_q$ and therefore the operation has an infinite number of solutions.
\item If $\pi$ is a plane parallel to the $xy$ coordinate plane, then $a=b=0$, $c\neq 0$, and $z'_q=-d/c$ (constant). At the same time, if conditions of item 1 above are not satisfied ($z'_q\neq -1$ or $P\neq Q$), then Eq.~(\ref{o63}) is quadratic in $y'_q$ and may have none to two solutions.
\item In any other case, Eq.~(\ref{o63}) results in a cubic for $y'_q$ or $z'_q$ and it may have one to three solutions.
\end{enumerate}
We conclude that the operation is well defined iff $P\neq Q$ or $m\not\subset \pi$ or $\pi$ is not perpendicular to the plane spanned by $P$ and $m$, and it may have none to three solutions.  An example is shown in Fig.~\ref{I5+I6}.

\subsection{Operation \texorpdfstring{$I_5+I_9$}{I5 + I9}}
\label{o7}
This operation combines two different constraints of codimension 2, and it may be stated as: \textit{Given lines $m$ and $n$, and a point $P$ not on $m$, fold along a plane to place $P$ onto  $m$ and to reflect half of $n$ onto its other half}. Its solution is a fold plane tangent to the parabolic cylinder generated by $P$ and $m$ (Fig.\ \ref{cylp}) and normal to $n$.

Again, assume $P$ and $m$ as in \S \ref{I5}, and $n$ a line with a direction vector $(a, b, c)$. A normal vector for the fold plane $\Delta$ is $(0, t, -2)$, and it must be parallel to the direction vector for $n$. Hence,
\begin{equation}
(0, t, -2) = k(a, b, c),
\label{o71}
\end{equation}
where $k$ is a coefficient. For a given line $m$ with $a=0$ and $c\neq 0$, the above equation produces $t=-2b/c$ which determines a unique fold plane. If $a\neq 0$ or $c=0$ then the equation does not have a solution.

\begin{figure}
\centering
\includegraphics{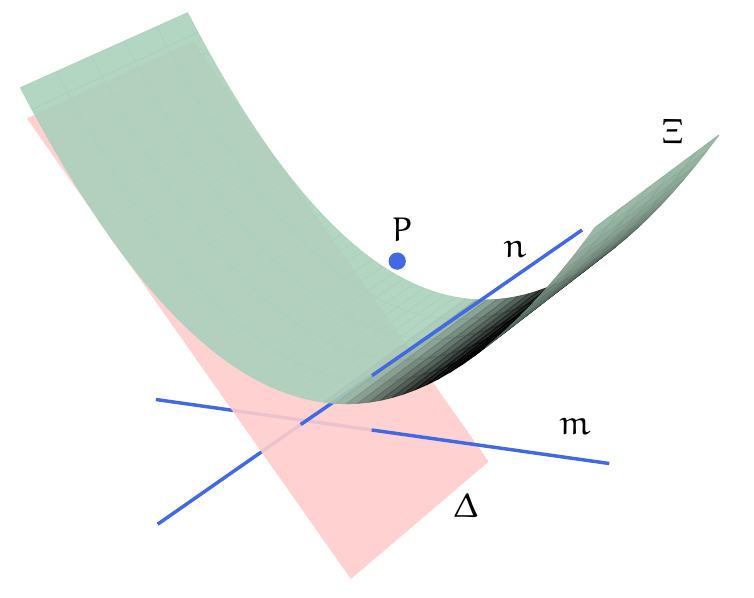}
\caption{Example of operation $I_5 + I_9$. The fold plane $\Delta$ is tangent to the parabolic cylinder $\Xi$, with focus $P$ and directrix line $m$, and normal to line $n$.}
\label{I5+I9}
\end{figure}

 We conclude that the operation has a solution iff line $n$ is parallel to the plane spanned by $P$ and $m$, and $n$ is not parallel to $m$, and that solution is unique. An example is shown in Fig. \ref{I5+I9}.

\subsection{Operation \texorpdfstring{$I_6 + I_8 + I_{11}$}{I6 + I8 + I11}}
This operation combines three different constraints of codimension 1, and it may be stated as: \textit{Given planes $\pi$ and $\tau$, a point $P$ not on $\pi$ and a point $Q$, fold along a plane passing through $Q$ to place $P$ onto  $\pi$ and to reflect half of $\tau$ onto its other half}. Its solution is a fold plane passing through $Q$, tangent to the paraboloid generated by $P$ and $\pi$ (Fig.\ \ref{paraboloid}), and normal to $\tau$.

Assume $P$ and $\pi$ as in \S \ref{I6}, point $Q$ at $(x_q, y_q, z_q)$ and $\tau$ described by $ax+by+cz+d=0$. The fold plane $\Delta$ is given by Eq. (\ref{deltaplane1}) and it must contain $Q$, therefore
\begin{equation}
2sx_q+2ty_q-4z_q = s^2+t^2.
\label{o6811}
\end{equation}
The normal vectors of $\Delta$, $(2s, 2t, -4)$, and $\tau$, $(a, b, c)$, must be mutually normal and so they satisfy
\begin{equation}
2sa+2tb-4c = 0.
\label{o6811b}
\end{equation}

If $a\neq 0$ or $b\neq0$ then the Eq.\ (\ref{o6811b}) may be solved for parameter $s$ or $t$, respectively, and substituting in Eq.\ (\ref{o6811}) produces a quadratic equation on the other parameter. Thus, the operation may have none to two solutions. If $a=0$ and $b=0$ (i.e., plane $\tau$ is parallel to $\pi$) then $c\neq 0$ and Eq.\  (\ref{o6811b}) is not satisfied.

\begin{figure}
\centering
\includegraphics{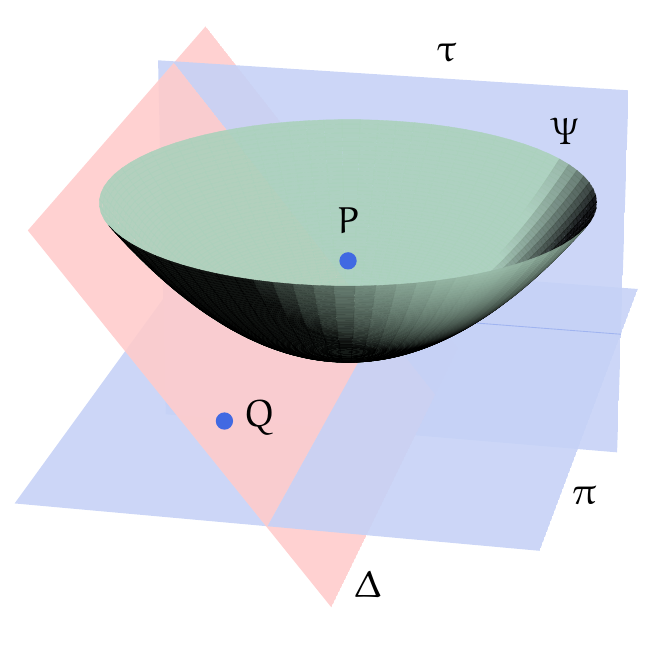}\hspace{1em}
\includegraphics{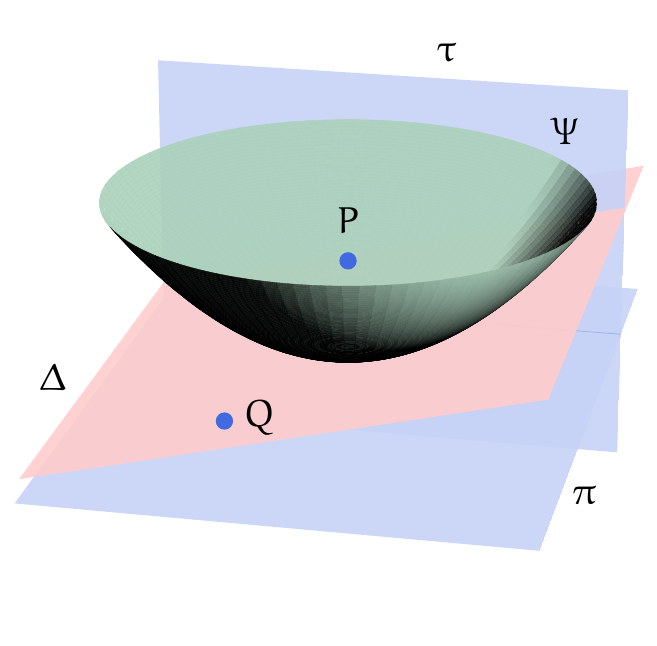}
\caption{Example of operation $I_6 + I_8 + I_{11}$ with two solutions. The fold plane $\Delta$ is tangent to paraboloid $\Psi$, with focus $P$ and directrix plane $\pi$, contains point $Q$ and is normal to plane $\tau$.}
\label{I6+I8+I11}
\end{figure}

An example is shown in Fig. \ref{I6+I8+I11}.

\subsection{Operation \texorpdfstring{$3I_6$}{3I6}}

This operation also combines three constraints of codimension 1, and it may be stated as: \textit{Given planes $\pi$, $\tau$, and $\rho$, a point $P$ not on $\pi$, a point $Q$ not on $\tau$, and a point $R$ not on $\rho$, fold along a plane to place $P$ onto  $\pi$, $Q$ onto  $\tau$, and $R$ onto  $\rho$}.
Its solution is a fold plane tangent to the three paraboloids defined by each of the incidences. 

Assume $P$ and $\pi$ as in \S \ref{I6}, point $Q$ at $(x_q, y_q, z_q)$ and $Q'$ at $(x'_q, y'_q, z'_q)$. The fold plane $\Delta$ is given by Eq. (\ref{deltaplane1}); it must be normal to segment $\overline{QQ'}$. Therefore, the normal vector of $\Delta$ is parallel to $\overline{QQ'}$ and satisfies
\begin{equation}
(s, t, -2) = k(x'_q-x_q, y'_q-y_q, z'_q-z_q),
\label{to61}
\end{equation}
where $k$ is a constant. Further, $\Delta$ passes through the midpoint of $\overline{QQ'}$ located at $((x_q +x'_q)/2, (y_q +y'_q)/2, (z_q+z'_q)/2)$. Replacing these coordinates into Eq. (\ref{deltaplane1}) produces
\begin{equation}
2s(x_q+x'_q)+ 2t(y_q+y'_q)-4(z_q+z'_q)=s^2+t^2.
\label{to62}
\end{equation}

Eliminating $s$, $t$, and $k$ from Eqs. (\ref{to61}) and (\ref{to62}) gives 
\begin{multline}
2(y_q-y'_q)^2 + 2(x_q-x'_q)^2=\\
-(x_q^2-{x'}_q^2)(z_q-z'_q)-(y_q^2-{y'}_q^2)(z_q-z'_q)-(z_q^2-{z'}_q^2)(z_q-z'_q).
\label{to63}
\end{multline}

Similarly, letting point $R$ be at $(x_r, y_r, z_r)$ and $R'$ at $(x'_r, y'_r, z'_r)$, produces
\begin{multline}
2(y_r-y'_r)^2 + 2(x_r-x'_r)^2=\\-(x_r^2-{x'}_r^2)(z_r-z'_r)-(y_r^2-{y'}_r^2)(z_r-z'_r)-(z_r^2-{z'}_r^2)(z_r-z'_r).
\label{to64}
\end{multline}

Segments $\overline{QQ'}$ and $\overline{RR'}$ are parallel (they are both normal to $\Delta$), therefore,

\begin{equation}
(x'_q-x_q, y'_q-y_q, z'_q-z_q) = \ell(x'_r-x_r, y'_r-y_r, z'_r-z_r).
\label{to64b}
\end{equation}
 
For planes $\tau$ and $\rho$ described by $a_\tau x+b_\tau y+c_\tau z+d_\tau=0$ and $a_\rho x+b_\rho y+c_\rho z+d_\rho=0$, respectively, where $a_\tau$, $b_\tau$, $c_\tau$, $d_\tau$, $a_\rho$, $b_\rho$, $c_\rho$, and $d_\rho$ are constants, the coordinates of $Q'$ and $R'$ satisfy
\begin{equation}
a_\tau x'_q+b_\tau y'_q+c_\tau z'_q+d_\tau =0,
\label{to65}
\end{equation}
and
\begin{equation}
a_\rho x'_q+b_\rho y'_q+c_\rho z'_q+d_\rho =0.
\label{to66}
\end{equation}

\begin{figure}
\centering
\includegraphics{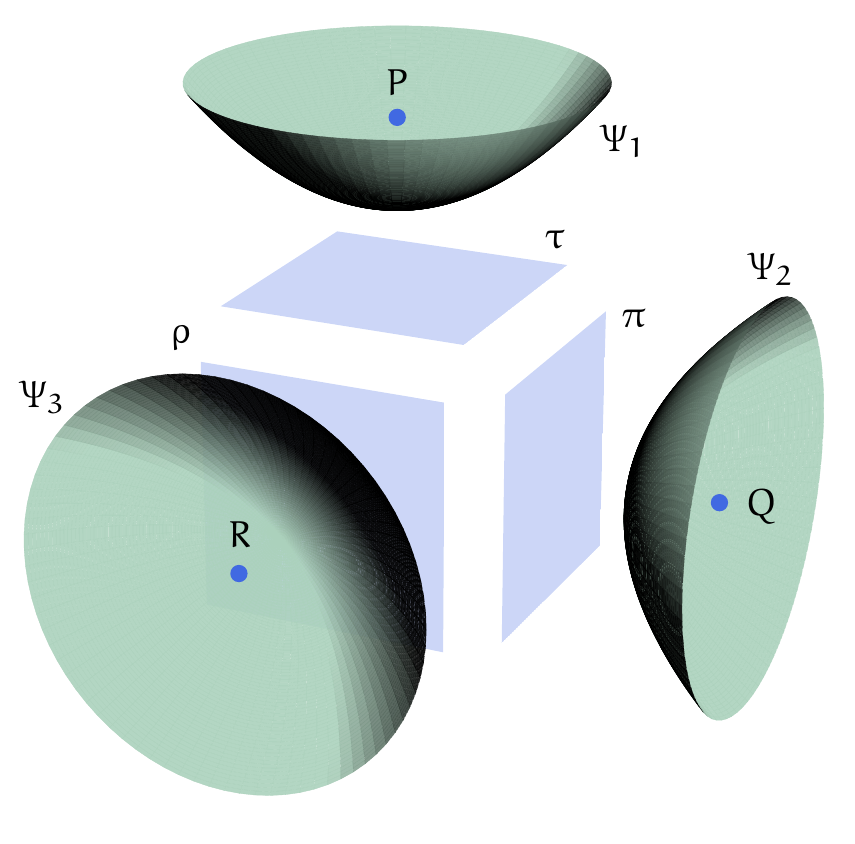}
\includegraphics{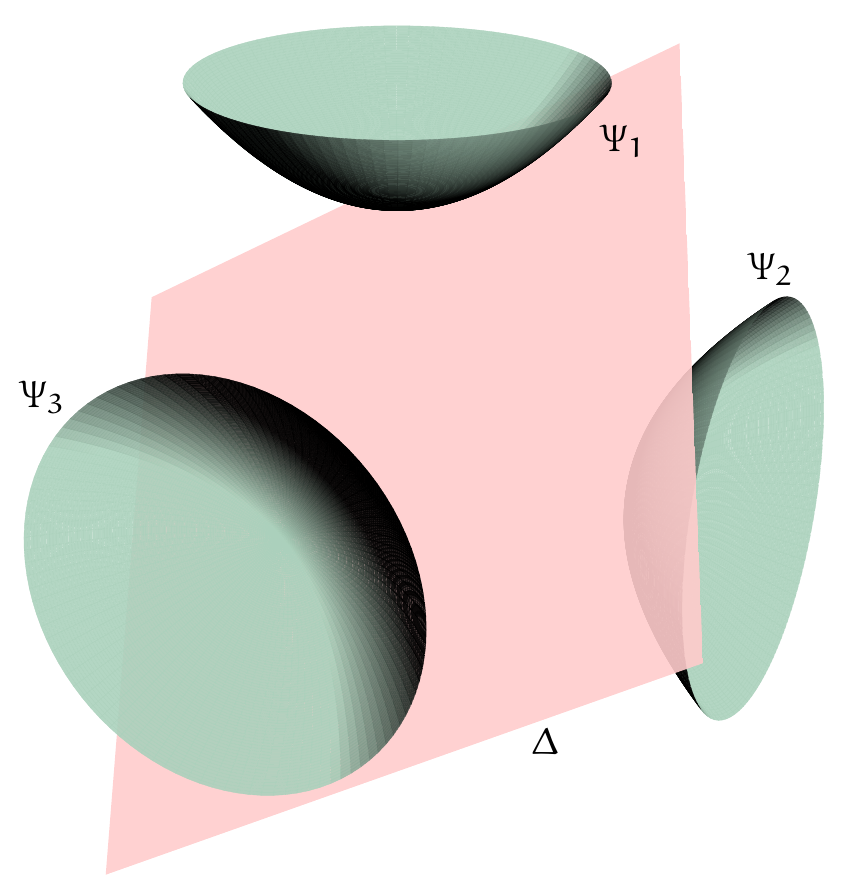}
\caption{Example of operation $3I_6$. Top: the fold planes are tangent to the three paraboloids $\Psi_1$, $\Psi_2$, and $\Psi_3$ with focuses and directrix planes $P$ and $\pi$, $Q$ and $\tau$, $R$ and $\rho$, respectively. Bottom: one of the seven solutions of the operation.}
\label{3I6}
\end{figure}

Equations (\ref{to63}) to (\ref{to66}) form a system of equations in $\ell$ and the coordinates of $Q'$ and $R'$. An equation for the fold plane $\Delta$ may be obtained by solving the system and next using the (real) solutions in Eqs.\ (\ref{to61}) and (\ref{deltaplane1}). According to Bezout's theorem \citep{Weisstein2018}, the number of solutions is at most the product of the equations' degrees, which results in an upper bound of 36. However, and analysis using Gr\"{o}bner basis computations \citep{sturmfels2002} in a Computer Algebra System reveals nine solutions at most. It may be also shown that two of the solutions are always complex, which leaves a maximum of seven possible solutions for the fold operation (calculations are lengthy and are therefore omitted here). An example is shown in Fig.~\ref{3I6}.

\section{Discussion}

A 3D extension of origami may be defined based on a set of at most 47 elementary folding operations (in principle, 50 operations may be defined, but three of them do not have a solution in an Euclidean space). For comparison, let us recall that 2D origami is based on seven elementary operations \citep{Alperin2006}  \citep[an additional eighth operation has also been considered by][]{Lucero2017}.  Folding takes place on a plane, and each operation satisfies a set of incidence constraints between given points, lines and planes with a finite number of solutions. Some of the operations have been analyzed, and it has been shown that one of them may have up to seven solutions. Again, for comparison, let us recall that the maximum number of solutions to 2D operations is 3. As a next step, a complete analysis of all 3D operations to determine conditions for solutions and their number would be desirable.

This analysis may found application to studies of Universe structure in cosmology. The concept of 3D origami tessellations has been used to model the large-scale distribution of matter in filaments, clusters and voids \citep{Neyrinck2014}. In that model, a 3D ``sheet'' of dark matter folds within a 6D Eulerian space (three spacial dimensions plus three velocity dimensions). Superposed regions of the folded 3D ``sheet'' form walls, streams and nodes, which approximates the observed distribution of galaxies. The present results might add to the mathematical formalization of such model and facilitate computer simulations.

\appendix
\section{Envelope of a family of surfaces}
\label{ApA}
The envelope of a one-parameter family of surfaces in space is a surface that touches each member of the family along a whole curve \citep{Courant1974,Hazewinkel1989}, and the set of curves of contact forms a one-parameter family that completely covers the envelope. 

Let the family be described by
 \begin{equation} 
F(t, x, y, z)=0,
\label{ApAe1}
\end{equation} 
where $t$ is a free parameter, $F\in C^1$ and $|F_x| + |F_y|+|F_z|\neq 0$. Then, a necessary condition for the existence on an envelope at a point $\mathbf{x}=(x, y, z)$ is that $\mathbf{x}$ satisfies both Eq.~(\ref{ApAe1}) and 
\begin{equation}
F_t(t, x, y, z)=0.
\label{ApAe2}
\end{equation}  
Eliminating the parameter $t$ from both equations, an equation called the discriminant is obtained, which must be satisfied by the envelope.
 
Sufficient conditions are $f\in C^2$ and, in addition to Eqs.~(\ref{ApAe1}
) and (\ref{ApAe2}), conditions
\begin{equation}
F_{tt}(t, x, y, z)\neq 0,
\label{ApAe3}
\end{equation}  
and
\begin{equation}
\left|\frac{D(F, F_t)}{D(x, y)}\right|+\left|\frac{D(F, F_t)}{D(y, z)}\right|+\left|\frac{D(F, F_t)}{D(z, x)}\right|\neq 0.
\end{equation}

In case of a two-parameter family of surfaces, the envelope touches each member of the family at one point.
If the family is described by
 \begin{equation} 
F(r, s, x, y, z)=0,
\label{ApAe4}
\end{equation} 
where $r$ and $s$ are free parameters, $F\in C^1$ and $|F_x| + |F_y|+|F_z|\neq 0$, then necessary conditions for the existence of the envelope at a point $(x, y, z)$ is that the point satisfies Eq.~(\ref{ApAe4}) and 
\begin{equation}
F_r(r, s, x, y, z)=0, \qquad F_s(r, s, x, y, z)=0
\label{ApAe5}
\end{equation}  
The discriminant equation is obtained by eliminating $r$ and $s$ from Eqs.~(\ref{ApAe4}) and (\ref{ApAe5}). 

Sufficient conditions are $f\in C^2$ and, in addition to Eqs.~(\ref{ApAe4}
) and (\ref{ApAe5}),
\begin{equation}
\frac{D(F, F_r, F_s)}{D(x, y,z)}\neq 0, \qquad \frac{D(F_r, F_s)}{D(r, s)}\neq 0. 
\label{ApAe6}
\end{equation}

\end{document}